# COMMUTING BIRTH-AND-DEATH PROCESSES

By Steven N. Evans[1], Bernd Sturmfels[2] and Caroline Uhler[3]

*University of California at Berkeley*

We use methods from combinatorics and algebraic statistics to study analogues of birth-and-death processes that have as their state space a finite subset of the $m$-dimensional lattice and for which the $m$ matrices that record the transition probabilities in each of the lattice directions commute pairwise. One reason such processes are of interest is that the transition matrix is straightforward to diagonalize, and hence it is easy to compute $n$ step transition probabilities. The set of commuting birth-and-death processes decomposes as a union of toric varieties, with the main component being the closure of all processes whose nearest neighbor transition probabilities are positive. We exhibit an explicit monomial parametrization for this main component, and we explore the boundary components using primary decomposition.

**1. Introduction.** Birth-and-death processes are among the simplest Markov chains. In discrete time, they model a particle that wanders back and forth on a sub-interval of the integers by taking unit size steps. Birth-and-death processes arise in fields ranging from ecology, where the points in the state space represent population size, and at each point in time only a single birth or death event can occur, to queuing theory where the states represent the number of individuals waiting in a queue and at any time one individual either joins or leaves the queue. The finite time behavior of a birth-and-death process is easy to study because the transition matrix $P$, whose entry $P(i, j)$ is the probability of going from state $i$ to state $j$ in one step, is tri-diagonal. This means that the matrix $P$ can be diagonalized using a related family of

Received December 2008.

[1]Supported in part by NSF Grants DMS-04-05778 and DMS-09-07630.

[2]Supported in part by NSF Grants DMS-04-56960 and DMS-07-57236.

[3]Supported by an International Fulbright Science and Technology Fellowship.

orthogonal polynomials, and this enables the computation of the power $P^n$ whose entries are the probabilities of going from one state to another in $n$ steps [13, 14].

It is natural to consider Markov chains that have as their state space products of intervals in higher-dimensional integer lattices. For instance, the ecology model could be extended to a situation in which individuals in the population have a type and one keeps track of the number of individuals of each type, or the queuing theory example generalizes to one where there are several servers, each with their own set of customers, and one follows the respective queue lengths.

There are also higher-dimensional models where one of the coordinates describes the quantity of primary interest while the others describe a fluctuating environment or background that modulates the dynamics of that quantity. This is the point of view taken in *quasi-birth-and-death* processes [15, 16], where the state of the primary variable is usually called the *level*, while the state of the subsidiary ones is the *phase*. Such models have been used, *inter alia*, to model queues in random environments and queues where the service times and inter-arrival times have phase-type distributions. The setting here is most often that of continuous time, but many of the same considerations apply to discrete time. A discussion of the discrete time case and its connection with matrix-valued orthogonal polynomials is given in [10]. There is also a huge literature in finance, economics and engineering on similar processes, where the terminology used is usually that of *regime switching* or *stochastic volatility* models. The setting there is, again, often in continuous time, and also the primary variable often has a continuous state space. However, numerical computations for such a model often involve approximation by one with discrete time and discrete state space. A few representative examples are [1, 4, 7, 11, 17, 19, 21, 22].

Unfortunately, even when such higher-dimensional models have only nearest neighbor transitions that are analogous to those of birth-and-death chains, it is no longer true that they are necessarily reversible, and so there is no hope that there will be a straightforward spectral decomposition like that afforded by the use of orthogonal polynomials. It is, therefore, natural to seek special cases where one can still recover something of the classical one-dimensional theory.

One case in which this is possible is if, when we write the transition matrix as the sum of matrices, one for the transitions in each coordinate direction, the resulting collection of matrices commute. We motivate our study by examining the two-dimensional case in Section 2. After deriving the algebraic constraints for commuting matrices, we show how spectral methods may still be used to compute $n$ step transition probabilities for a process on a two-dimensional grid. The same approach extends without difficulty to higher-dimensional grids.

With this observation in mind, we aim to understand what restrictions are placed on a Markov chain by the requirement that the matrices appearing in such a decomposition commute. In Section 3, we obtain an essentially unique parametrization of such a commuting model under the assumption that all nearest neighbor transition probabilities are positive; and, under the same assumption, we characterize the minimal number of constraints on the transition probabilities that are necessary to ensure commutation. We then use methods from linear algebra, matroid theory and algebraic statistics to further analyze the binomial ideal generated by the commutation condition. In Section 4, we determine for which grid graphs the matrix describing the parametrization is unimodular. For these lattices, we derive an explicit graphical representation of the Graver basis, which is the most inclusive in the hierarchy of lattice bases in [5], Section 1.3. In Section 5, we explore the extraneous boundary components of our commuting variety. These correspond to families of commuting birth-and-death processes having some zero transition probabilities which are not limits of commuting birth-and-death processes with positive transition probabilities. This underlines the applicability of binomial primary decomposition in probability and statistics, well beyond the contexts envisioned by Diaconis, Eisenbud and Sturmfels in [3].

## 2. Motivation: Birth-and-death processes in dimension two.

Consider a discrete-time, time-homogeneous Markov chain $Z = (Z_k)_{k=0}^{\infty}$ that has as its state space the finite two-dimensional grid $E := \{0, 1, \ldots, m\} \times \{0, 1, \ldots, n\}$. Suppose the chain makes jumps of size one either upwards, downwards, to the right or to the left. In other words, if we impose a graph structure on $E$ by declaring that two states $(i', j')$ and $(i'', j'')$ are connected by an edge if and only if $|i' - i''| + |j' - j''| = 1$, then the chain can only make "nearest neighbor" jumps. If we draw $E$ in the plane and include the edges as intervals of unit length, then the resulting figure is made up of $m \times n$ squares, and so we refer to $E$ as the $m \times n$ *grid*.

The dynamics of $Z$ are specified by the transition probabilities

$$L_{i,j} := \mathbb{P}\{Z_{k+1} = (i-1, j) \mid Z_k = (i, j)\},$$

$$R_{i,j} := \mathbb{P}\{Z_{k+1} = (i+1, j) \mid Z_k = (i, j)\},$$

$$D_{i,j} := \mathbb{P}\{Z_{k+1} = (i, j-1) \mid Z_k = (i, j)\},$$

$$U_{i,j} := \mathbb{P}\{Z_{k+1} = (i, j+1) \mid Z_k = (i, j)\}.$$

These transition probabilities are nonnegative real numbers that satisfy $L_{i,j} + R_{i,j} + D_{i,j} + U_{i,j} \leq 1$. This inequality is allowed to be strict. As usual, strict inequality is interpreted in terms of an adjoined absorbing state $\dagger$ with $\mathbb{P}\{Z_{k+1} = \dagger \mid Z_k = (i, j)\} = 1 - (L_{i,j} + R_{i,j} + D_{i,j} + U_{i,j})$.

The transition matrix $P$ of the Markov chain $Z$ has format $[(m+1)(n+1)] \times [(m+1)(n+1)]$. It may be written as a sum of two matrices, $P =$ $P_h + P_v$, one for the horizontal moves and one for the vertical moves. The horizontal matrix $P_h$ commutes with the vertical matrix $P_v$ if and only if the following four constraints hold for all index pairs $(i, j) \in \{0, 1, \ldots, m - 1\} \times \{0, 1, \ldots, n - 1\}$:

$$
\begin{aligned}
U_{i,j}R_{i,j+1} &= R_{i,j}U_{i+1,j} &&\text{(up-right);} \\
D_{i,j+1}R_{i,j} &= R_{i,j+1}D_{i+1,j+1} &&\text{(down-right);} \\
D_{i+1,j+1}L_{i+1,j} &= L_{i+1,j+1}D_{i,j+1} &&\text{(down-left);} \\
U_{i+1,j}L_{i+1,j+1} &= L_{i+1,j}U_{i,j} &&\text{(up-left).}
\end{aligned}
$$

(2.1)

In English, for each corner of a square in the grid, the probability of going from that corner to the diagonally opposite corner of the square in two steps is the same for the two possible paths.

EXAMPLE 2.1 ($m = 2$ and $n = 1$). We label the vertices of the $2 \times 1$ grid as follows.

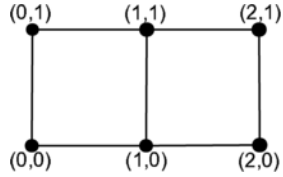

The $6 \times 6$ transition matrix $P = P_h + P_v$ is the sum of the horizontal transition matrix

$$
P_h = \begin{array}{c|cccccc}
 & (0,0) & (0,1) & (1,0) & (1,1) & (2,0) & (2,1) \\
\hline
(0,0) & 0 & 0 & R_{00} & 0 & 0 & 0 \\
(0,1) & 0 & 0 & 0 & R_{01} & 0 & 0 \\
(1,0) & L_{10} & 0 & 0 & 0 & R_{10} & 0 \\
(1,1) & 0 & L_{11} & 0 & 0 & 0 & R_{11} \\
(2,0) & 0 & 0 & L_{20} & 0 & 0 & 0 \\
(2,1) & 0 & 0 & 0 & L_{21} & 0 & 0
\end{array}
$$

and the vertical transition matrix

$$
P_v = \begin{array}{c|cccccc}
 & (0,0) & (0,1) & (1,0) & (1,1) & (2,0) & (2,1) \\
\hline
(0,0) & 0 & U_{00} & 0 & 0 & 0 & 0 \\
(0,1) & D_{01} & 0 & 0 & 0 & 0 & 0 \\
(1,0) & 0 & 0 & 0 & U_{10} & 0 & 0 \\
(1,1) & 0 & 0 & D_{11} & 0 & 0 & 0 \\
(2,0) & 0 & 0 & 0 & 0 & 0 & U_{20} \\
(2,1) & 0 & 0 & 0 & 0 & D_{21} & 0
\end{array}.
$$

Suppose that the 14 transition probabilities $L_{ij}, R_{ij}, D_{ij}, U_{ij}$ are all strictly positive. Then the commuting relations (2.1) are found to be equivalent to the following rank constraints:

(2.2)
$$\operatorname{rank} \begin{pmatrix} R_{00} & U_{00} & L_{11} & D_{11} \\ R_{01} & U_{10} & L_{10} & D_{01} \end{pmatrix} = 1 \quad \text{and}$$

$$\operatorname{rank} \begin{pmatrix} R_{10} & U_{10} & L_{21} & D_{21} \\ R_{11} & U_{20} & L_{20} & D_{11} \end{pmatrix} = 1.$$

Indeed, the eight equations in (2.1) are among the twelve $2 \times 2$-minors of these two $2 \times 4$-matrices. The constraints (2.2) imply that our commuting variety $\{P_h P_v = P_v P_h\}$ has the parametrization

(2.3)
$$R_{00} = h_1 \frac{a_{00}}{a_{10}}, \qquad R_{01} = h_1 \frac{a_{01}}{a_{11}}, \qquad R_{10} = h_2 \frac{a_{10}}{a_{20}}, \qquad R_{11} = h_2 \frac{a_{11}}{a_{21}},$$

$$L_{10} = h_1 \frac{a_{10}}{a_{00}}, \qquad L_{11} = h_1 \frac{a_{11}}{a_{01}}, \qquad L_{20} = h_2 \frac{a_{20}}{a_{10}}, \qquad L_{21} = h_2 \frac{a_{21}}{a_{11}},$$

$$U_{00} = v_1 \frac{a_{00}}{a_{01}}, \qquad U_{10} = v_1 \frac{a_{10}}{a_{11}}, \qquad U_{20} = v_1 \frac{a_{20}}{a_{21}},$$

$$D_{01} = v_1 \frac{a_{01}}{a_{00}}, \qquad D_{11} = v_1 \frac{a_{11}}{a_{10}}, \qquad D_{21} = v_1 \frac{a_{21}}{a_{20}}.$$

An analogous parametrization for commuting birth-and-death processes on larger two-dimensional grids and in higher dimensions will be derived in Section 3 and further studied in Section 4. The analogues of the parameters $h_1, h_2, v_1$ and the parameters $a_{i,j}$ are straightforward to compute given the transition matrix $P$. We note that if the transition probabilities are allowed to be zero then (2.3) is insufficient because (2.1) can be satisfied even if (2.2) fails. The identification of such boundary phenomena is the content of our algebraic discussion in Section 5.

The parametrization (2.3) of the transition matrix $P$ in Example 2.1 can be rewritten as

(2.4)
$$P_h = A Q_h A^{-1} \quad \text{and} \quad P_v = A Q_v A^{-1},$$

where the three matrices appearing on the right-hand side are

$$A = \begin{array}{c|cccccc} & (0,0) & (0,1) & (1,0) & (1,1) & (2,0) & (2,1) \\ \hline (0,0) & a_{00} & 0 & 0 & 0 & 0 & 0 \\ (0,1) & 0 & a_{01} & 0 & 0 & 0 & 0 \\ (1,0) & 0 & 0 & a_{10} & 0 & 0 & 0 \\ (1,1) & 0 & 0 & 0 & a_{11} & 0 & 0 \\ (2,0) & 0 & 0 & 0 & 0 & a_{20} & 0 \\ (2,1) & 0 & 0 & 0 & 0 & 0 & a_{21} \end{array},$$

$$Q_h = \begin{array}{c|cccccc} & (0,0) & (0,1) & (1,0) & (1,1) & (2,0) & (2,1) \\ \hline (0,0) & 0 & 0 & h_1 & 0 & 0 & 0 \\ (0,1) & 0 & 0 & 0 & h_1 & 0 & 0 \\ (1,0) & h_1 & 0 & 0 & 0 & h_2 & 0 \\ (1,1) & 0 & h_1 & 0 & 0 & 0 & h_2 \\ (2,0) & 0 & 0 & h_2 & 0 & 0 & 0 \\ (2,1) & 0 & 0 & 0 & h_2 & 0 & 0 \end{array},$$

$$Q_v = \begin{array}{c|cccccc} & (0,0) & (0,1) & (1,0) & (1,1) & (2,0) & (2,1) \\ \hline (0,0) & 0 & v_1 & 0 & 0 & 0 & 0 \\ (0,1) & v_1 & 0 & 0 & 0 & 0 & 0 \\ (1,0) & 0 & 0 & 0 & v_1 & 0 & 0 \\ (1,1) & 0 & 0 & v_1 & 0 & 0 & 0 \\ (2,0) & 0 & 0 & 0 & 0 & 0 & v_1 \\ (2,1) & 0 & 0 & 0 & 0 & v_1 & 0 \end{array}.$$

As we shall see in Theorem 3.1, the matrix parametrization (2.4) is valid for the $m \times n$-grid. Here $A$ is a diagonal matrix of size $(m+1)(n+1)$, the matrix $Q_v$ is block diagonal with $m+1$ identical $(n+1) \times (n+1)$-blocks $R_v$ that are symmetric and tri-diagonal and $Q_h$ is block diagonal with $n+1$ identical $(m+1) \times (m+1)$-blocks $R_h$ that are symmetric and tri-diagonal. Note that in order to make the block diagonal structure of the matrix $Q_h$ apparent, it is necessary to re-order the rows and columns as follows: $(0,0), (1,0), (2,0), (0,1), (1,1), (2,1)$. Thus the tri-diagonal matrices $R_h$ and $R_v$ satisfy

$$Q_h((i', j'), (i'', j'')) = \begin{cases} R_h(i', i''), & \text{if } j' = j'', \\ 0, & \text{otherwise,} \end{cases}$$

and

$$Q_v((i', j'), (i'', j'')) = \begin{cases} R_v(j', j''), & \text{if } i' = i'', \\ 0, & \text{otherwise.} \end{cases}$$

By the spectral theorem for real symmetric matrices [12], Theorem 4.1.5, the $r$th power of the matrix $R_h$ has entries

$$R_h^r(i', i'') = \sum_{k=0}^{m} \lambda_{h,k}^r u_{h,k}(i') u_{h,k}(i''),$$

where $\lambda_{h,0}, \lambda_{h,1}, \ldots, \lambda_{h,m}$ are the eigenvalues of $R_h$ and the $u_{h,k}$ are the corresponding orthonormalized eigenvectors. With analogous notation, for any positive integer $r$, we have

$$R_v^r(i', i'') = \sum_{\ell=0}^{n} \lambda_{v,\ell}^r u_{v,\ell}(j') u_{v,\ell}(j'').$$

Now, because $P_h$ and $P_v$ commute, the $t$-step transition probability matrix $P^t$ is given by

$$P^t = (P_h + P_v)^t = \sum_{s=0}^{t} \binom{t}{s} P_h^s P_v^{t-s}$$

$$= \sum_{s=0}^{t} \binom{t}{s} (A Q_h A^{-1})^s (A Q_v A^{-1})^{t-s}$$

$$= A \left[ \sum_{s=0}^{t} \binom{t}{s} Q_h^s Q_v^{t-s} \right] A^{-1}.$$

By combining the above formulas, we obtain the following simple expression for the $t$-step transition probabilities of the Markov chain $Z$:

$$P^t((i', j'), (i'', j''))$$
$$= a_{i',j'} \left[ \sum_k \sum_\ell (\lambda_{h,k} + \lambda_{v,\ell})^t u_{h,k}(i') u_{h,k}(i'') u_{v,\ell}(j') u_{v,\ell}(j'') \right] a_{i'',j''}^{-1}.$$

Moreover, note that for a constant $c > 0$ and $I$ the identity matrix,

$$(cI + P)^t((i', j'), (i'', j''))$$
$$= a_{i',j'} \left[ \sum_k \sum_\ell (c + \lambda_{h,k} + \lambda_{v,\ell})^t u_{h,k}(i') u_{h,k}(i'') u_{v,\ell}(j') u_{v,\ell}(j'') \right] a_{i'',j''}^{-1},$$

so we may readily compute the $t$-step transition probabilities of certain chains for which the particle has a constant positive probability of not moving from its current position at each step. The calculations we have just done extend to the finite higher-dimensional grids $E$ to be discussed in the next section. The problem of finding the $t$-step transition probabilities for these models reduces to finding the eigenvalues and eigenvectors of a collection of symmetric tri-diagonal matrices whose nonzero entries are the parameters $W(u, v)$ in Theorem 3.1.

In our definition of the Markov chain $Z$ above, we did not require that the rows of the transition matrix $P$ add to 1. In order to get a feeling for the effect of imposing that extra constraint, we return to the case $m = 2$ and $n = 1$ in Example 2.1. The rows of $P$ add to 1 if and only if 1 is an eigenvalue of the symmetric matrix

$$Q_h + Q_v =$$

|       | $(0,0)$ | $(0,1)$ | $(1,0)$ | $(1,1)$ | $(2,0)$ | $(2,1)$ |
|-------|---------|---------|---------|---------|---------|---------|
| $(0,0)$ | $0$   | $v_1$   | $h_1$   | $0$     | $0$     | $0$     |
| $(0,1)$ | $v_1$ | $0$     | $0$     | $h_1$   | $0$     | $0$     |
| $(1,0)$ | $h_1$ | $0$     | $0$     | $v_1$   | $h_2$   | $0$     |
| $(1,1)$ | $0$   | $h_1$   | $v_1$   | $0$     | $0$     | $h_2$   |
| $(2,0)$ | $0$   | $0$     | $h_2$   | $0$     | $0$     | $v_1$   |
| $(2,1)$ | $0$   | $0$     | $0$     | $h_2$   | $v_1$   | $0$     |

with corresponding right, and hence also left, eigenvector $(a_{00}^{-1}, a_{01}^{-1}, a_{10}^{-1}, a_{11}^{-1}, a_{20}^{-1}, a_{21}^{-1})$.

The preceding matrix is irreducible when the parameters $h_1, h_2, v_1$ are positive, and the Perron–Frobenius theorem [12], Theorem 8.4.4, guarantees that the eigenvalue with largest absolute value is positive, with an eigenvector that has positive entries and is unique up to a constant multiple. Consequently, if positive parameters $h_1, h_2, v_1$ are given, then replacing them by $bh_1, bh_2, bv_1$ for a suitable positive constant $b$ permits the choice of parameters $a_{ij}$ that are unique up to an irrelevant constant multiple such that the rows of $P$ sum to 1. More generally, replacing $h_1, h_2, v_1$ by $bh_1, bh_2, bv_1$ for a suitable positive constant $b$ permits the choice of parameters $a_{ij}$ that are unique up to an irrelevant constant multiple such that the rows of $P$ sum to $1 - c$ for a given constant $0 < c < 1$, in which case the rows of the matrix $cI + P$ sum to 1. These considerations extend to larger two-dimensional grids and to higher dimensions.

## 3. Parametrization.

Suppose now that we have a discrete-time, time-homogeneous Markov chain with state space

$$E := \{0, 1, \ldots, n_1\} \times \{0, 1, \ldots, n_2\} \times \cdots \times \{0, 1, \ldots, n_m\}.$$

The $m$-dimensional grid $E$ indexes the rows and columns of the transition matrix $P$. We assume that $P(u, v) = 0$ unless $u \sim v$ by which we mean that $u - v \in \{\pm e_k\}$ for some $k \in \{1, \ldots, m\}$. Here $e_k$ is the standard basis vector with entry 1 in the $k$th coordinate and 0 elsewhere. Thus $P$ describes a Markov chain with nearest neighbor transitions on the graph with vertices $E$ where two vertices $u$ and $v$ are connected by an edge if $u \sim v$. For $1 \le k \le m$, define a matrix $P_k$ by

$$P_k(u, v) := \begin{cases} P(u, v), & \text{if } u - v \in \{\pm e_k\}, \\ 0, & \text{otherwise.} \end{cases}$$

These matrices are the analogues of the matrices $P_h$ and $P_v$ for the two-dimensional grid. The requirement that the matrices $P_i$ and $P_j$ commute for all pairs of indices $1 \leq i < j \leq m$ is equivalent to the condition that the following quadratic expressions vanish for all $i, j$ and $u$:

$$
\begin{aligned}
& P(u, u + e_i) P(u + e_i, u + e_i + e_j) \\
& \quad - P(u, u + e_j) P(u + e_j, u + e_i + e_j), \\
& P(u, u + e_i) P(u + e_i, u + e_i - e_j) \\
& \quad - P(u, u - e_j) P(u - e_j, u + e_i - e_j), \\
& P(u, u - e_i) P(u - e_i, u - e_i + e_j) \\
& \quad - P(u, u + e_j) P(u + e_j, u - e_i + e_j), \\
& P(u, u - e_i) P(u - e_i, u - e_i - e_j) \\
& \quad - P(u, u - e_j) P(u - e_j, u - e_i - e_j).
\end{aligned}
$$

(3.1)

Our aim is to solve this system of polynomial equations for the unknowns $P(u, v)$. The next theorem offers such a solution under the assumption that the unknowns are strictly positive.

THEOREM 3.1. *Suppose that $P(u, v) > 0$ for all $u, v \in E$ such that $u \sim v$. Then the matrices $P_1, P_2, \ldots, P_m$ commute pairwise if and only if*

$$
(3.2) \qquad P(u, v) = a_u \cdot W(u, v) \cdot a_v^{-1}, \qquad u, v \in E,
$$

*for some collection of constants $a_u$ and $W(u, v)$, $u, v \in E$, that satisfy $W(u', v') = W(u'', v'')$ when $v' - u' \in \{\pm e_k\}$ for some $k$ and $\{u'', v''\} = \{u' + w, v' + w\}$ for some $w \in \sum_{\ell \neq k} \mathbb{Z} e_\ell$. The constants $a_u$ are unique up to a common multiple, and the constants $W(u, v)$ are unique.*

EXAMPLE 3.2 ($2 \times 1$ grid). If $m = 2$, $n_1 = 2$ and $n_2 = 1$, then this is the parametrization (2.3) in Example 2.1 with $h_1 = W((0,0), (1,0))$, $h_2 = W((1,0), (2,0))$ and $v_1 = W((0,0), (0,1))$.

PROOF OF THEOREM 3.1. The sufficiency of the stated condition is straightforward. For example,

$$
\begin{aligned}
& P(u, u + e_i) P(u + e_i, u + e_i + e_j) \\
& \quad = a_u W(u, u + e_i) a_{u+e_i}^{-1} a_{u+e_i} W(u + e_i, u + e_i + e_j) a_{u+e_i+e_j}^{-1} \\
& \quad = a_u W(u, u + e_j) a_{u+e_j}^{-1} a_{u+e_j} W(u + e_j, u + e_i + e_j) a_{u+e_i+e_j}^{-1} \\
& \quad = P(u, u + e_j) P(u + e_j, u + e_i + e_j),
\end{aligned}
$$

because $W(u, u + e_i) = W(u + e_j, u + e_i + e_j)$ and $W(u, u + e_j) = W(u + e_i, u + e_i + e_j)$.

For the converse, we first show that if the matrices $P_1, \ldots, P_m$ commute, then the transition matrix $P$ is reversible; that is, there are positive constants $b_u$ such that $b_u P(u, v) = b_v P(v, u)$ for all $u, v \in E$. Write $\|u\| = u_1 + \cdots + u_m$ for $u = (u_1, \ldots, u_m) \in E$. We claim it is possible to construct the constants $b_u$, $u \in E$, by induction on the value of $\|u\|$, and that their values are determined once $b_0$ is specified. For $u \in E$ with $\|u\| = 1$, we must take $b_u = b_0 P(0, u)/P(u, 0)$. For $v \in E$ with $\|v\| = 2$, we need to be able to find $b_v$ such that $b_u P(u, v) = b_v P(v, u)$ for all $u \in E$ with $\|u\| \leq 1$. This is equivalent to showing that we can choose $b_v$ so that $b_u P(u, v) = b_v P(v, u)$ for all $u \in E$ of the form $u = v - e_k$ for some $k$. If $v = 2e_k$ for some $k$, then we must take

$$b_v = b_{e_k} \frac{P(e_k, 2e_{2k})}{P(2e_k, e_k)} = b_0 \frac{P(0, e_k)}{P(e_k, 0)} \frac{P(e_k, 2e_{2k})}{P(2e_k, e_k)}.$$

If $v = e_i + e_j$ for $i \neq j$, then it is necessary that

$$b_{e_i} \frac{P(e_i, e_i + e_j)}{P(e_i + e_j, e_i)} = b_0 \frac{P(0, e_i)}{P(e_i, 0)} \frac{P(e_i, e_i + e_j)}{P(e_i + e_j, e_i)}$$

and

$$b_{e_j} \frac{P(e_j, e_i + e_j)}{P(e_i + e_j, e_j)} = b_0 \frac{P(0, e_j)}{P(e_j, 0)} \frac{P(e_j, e_i + e_j)}{P(e_i + e_j, e_j)}$$

are equal, in which case we must take $b_v$ to be the common value. However, from (3.1) we have

$$P(0, e_i) P(e_i, e_i + e_j) = P(0, e_j) P(e_j, e_i + e_j)$$

and

$$P(e_i + e_j, e_i) P(e_i, 0) = P(e_i + e_j, e_j) P(e_j, 0).$$

Hence the value of $b_v$ is uniquely defined. Continuing in this way shows that if $w \in E$ and $0 = w_0, w_1, \ldots, w_N = w$ with $w_{M+1} - w_M \in \{e_1, \ldots, e_m\}$ for $0 \leq M \leq N - 1$, then the value of

$$\prod_{M=0}^{N-1} \frac{P(w_M, w_{M+1})}{P(w_{M+1}, w_M)}$$

does not depend on the choice of $w_1, \ldots, w_{N-1}$. Moreover, if we take $b_w$ to be this common value, then the collection of constants $b_u$, $u \in E$, is such that $b_u P(u, v) = b_v P(v, u)$ for all $u, v \in E$, and this is the unique collection with that property and the given value of $b_0$.

Now, suppose that we have positive constants $b_u$, $u \in E$, such that $b_u P(u, v) = b_v P(v, u)$ for all $u, v \in E$. If we set $a_u = b_u^{-1/2}$ and define $W(u, v)$ by

$$(3.3) \qquad P(u, v) = a_u W(u, v) a_v^{-1},$$

then

(3.4)  $$W(u,v) = W(v,u).$$

Conversely, if we have constants $a_u$ and $W(u,v)$ that satisfy (3.3) and (3.4), then

$$\frac{P(u,v)}{P(v,u)} = \frac{a_u W(u,v) a_v^{-1}}{a_v W(v,u) a_u^{-1}} = a_u^2 a_v^{-2},$$

when $u \sim v$, and so $a_u^{-2} P(u,v) = a_v^{-2} P(v,u)$, $u,v \in E$. From what we have argued above, this implies that there are constants, $a_u$, $u \in E$ and $W(u,v)$, $u,v \in E$, that satisfy (3.3) and (3.4); the $a_u$ are unique up to a common constant multiple, and the $W(u,v)$ are unique.

To complete the proof, we need to check that the $W(u,v)$ have the additional properties listed in the statement of the theorem. Because $P(u,v)$, $u,v \in E$ is a common zero of the polynomials in (3.1), it follows that $W(u,v)$, $u,v \in E$ is a common zero of the same set of polynomials. The constraints that are associated with a particular two-dimensional face of one of the hypercubes in the grid with vertices $\{u, u+e_i, u+e_j, u+e_i+e_j\}$ are

$$W(u, u+e_i)W(u+e_i, u+e_i+e_j)$$
$$= W(u, u+e_j)W(u+e_j, u+e_i+e_j),$$

$$W(u+e_i, u+e_i+e_j)W(u+e_i+e_j, u+e_j)$$
$$= W(u+e_i, u)W(u, u+e_j),$$

$$W(u+e_i+e_j, u+e_j)W(u+e_j, u)$$
$$= W(u+e_i+e_j, u+e_i)W(u+e_i, u),$$

$$W(u+e_j, u)W(u, u+e_i)$$
$$= W(u+e_j, u+e_i+e_j)W(u+e_i+e_j, u+e_i).$$

Because $W(u,v) = W(v,u)$ for all $u,v \in E$, these four constraints are equivalent to

$$W(u, u+e_i)W(u+e_i, u+e_i+e_j) = W(u, u+e_j)W(u+e_j, u+e_i+e_j)$$

and

$$W(u+e_i, u+e_i+e_j)W(u+e_j, u+e_i+e_j) = W(u, u+e_i)W(u, u+e_j)$$

and hence to the equations

$$W(u, u+e_i) = W(u+e_j, u+e_i+e_j)$$

and

$$W(u, u+e_j) = W(u+e_i, u+e_j+e_i).$$

Iterating these two equations yields the remaining assertions stated in Theorem 3.1. □

REMARK 3.3. Note that if $u - v \in \{\pm e_k\}$, then $\{u, v\} = \{he_k + r, (h + 1)e_k + r\}$ for unique values of $h \in \{0, 1, \ldots, n_k - 1\}$ and $r \in \sum_{\ell \neq k} \mathbb{Z}e_\ell$. Therefore, if the necessary and sufficient conditions of Theorem 3.1 hold, then $W(u, v) = W(he_k, (h + 1)e_k)$. Thus the parametrization involves $n_1 + n_2 + \cdots + n_m$ uniquely defined parameters of this form and $(n_1 + 1)(n_2 + 1) \cdots (n_m + 1)$ parameters of the form $a_u$. The latter are uniquely defined up to a common multiple.

The set of constraints (3.1) has some redundancies when the unknown quantities $P(u, v)$, $u \sim v$, are positive. For example, any three of the equations in (2.1) implies the fourth. In what follows, we present a linear algebra approach to identifying and representing constraints that are independent. Note first that the quadratic equations in (3.1) are all of the form $P(a, b)P(b, c) - P(a, d)P(d, c) = 0$. This condition is equivalent to $[P(a, b)P(b, c)]/[P(a, d)P(d, c)] = 1$ or $Q(a, b) + Q(b, c) - Q(a, d) - Q(d, c) = 0$ where we write $Q(u, v)$ for the logarithm of $P(u, v)$. This suggests the following encoding of (3.1) as rows of a matrix with entries from $\{-1, 0, +1\}$.

NOTATION 3.4. Let $\mathcal{S}^{(n_1, \ldots, n_m)}$ be the matrix that has one row for each polynomial in (3.1) and columns indexed by ordered pairs $(u, v) \in E \times E$ with $u \sim v$, with the row corresponding to a polynomial of the form $P(a, b)P(b, c) - P(a, d)P(d, c)$ having entries $+1, +1, -1, -1$ in the columns associated with the pairs $(a, b), (b, c), (a, d), (d, c)$ and zero entries elsewhere.

REMARK 3.5. The matrix $\mathcal{S}^{(n_1, \ldots, n_m)}$ has format

$$\left( 4 \sum_{1 \leq i < j \leq m} \left( n_i n_j \prod_{k \neq i, j} (n_k + 1) \right) \right) \times \left( 2 \sum_{i=1}^{m} n_i \prod_{j \neq i} (n_j + 1) \right).$$

We wish to determine the rank of $\mathcal{S}^{(n_1, \ldots, n_m)}$, that is, the dimension of the vector space spanned by the rows, as this gives the size of a maximal set of independent constraints for positive $P(u, v)$.

EXAMPLE 3.6 ($2 \times 1$ grid constraint matrix). The $8 \times 14$-matrix $\mathcal{S}^{(2,1)}$ has rank 6 and equals

|  | $R_{00}$ | $R_{01}$ | $R_{10}$ | $R_{11}$ | $L_{10}$ | $L_{11}$ | $L_{20}$ | $L_{21}$ |
|---|---|---|---|---|---|---|---|---|
| $R_{00}U_{10} - U_{00}R_{01}$ | 1 | $-1$ | 0 | 0 | 0 | 0 | 0 | 0 |
| $U_{10}L_{11} - L_{10}U_{00}$ | 0 | 0 | 0 | 0 | $-1$ | 1 | 0 | 0 |
| $L_{11}D_{01} - D_{11}L_{10}$ | 0 | 0 | 0 | 0 | $-1$ | 1 | 0 | 0 |
| $D_{01}R_{00} - R_{01}D_{11}$ | 1 | $-1$ | 0 | 0 | 0 | 0 | 0 | 0 |
| $R_{10}U_{20} - U_{10}R_{11}$ | 0 | 0 | 1 | $-1$ | 0 | 0 | 0 | 0 |
| $U_{20}L_{21} - L_{20}U_{10}$ | 0 | 0 | 0 | 0 | 0 | 0 | $-1$ | 1 |
| $L_{21}D_{11} - D_{21}L_{20}$ | 0 | 0 | 0 | 0 | 0 | 0 | $-1$ | 1 |
| $D_{11}R_{10} - R_{11}D_{21}$ | 0 | 0 | 1 | $-1$ | 0 | 0 | 0 | 0 |

|  | $U_{00}$ | $U_{10}$ | $U_{20}$ | $D_{01}$ | $D_{11}$ | $D_{21}$ |
|---|---|---|---|---|---|---|
| $R_{00}U_{10} - U_{00}R_{01}$ | $-1$ | 1 | 0 | 0 | 0 | 0 |
| $U_{10}L_{11} - L_{10}U_{00}$ | $-1$ | 1 | 0 | 0 | 0 | 0 |
| $L_{11}D_{01} - D_{11}L_{10}$ | 0 | 0 | 0 | 1 | $-1$ | 0 |
| $D_{01}R_{00} - R_{01}D_{11}$ | 0 | 0 | 0 | 1 | $-1$ | 0 |
| $R_{10}U_{20} - U_{10}R_{11}$ | 0 | $-1$ | 1 | 0 | 0 | 0 |
| $U_{20}L_{21} - L_{20}U_{10}$ | 0 | $-1$ | 1 | 0 | 0 | 0 |
| $L_{21}D_{11} - D_{21}L_{20}$ | 0 | 0 | 0 | 0 | 1 | $-1$ |
| $D_{11}R_{10} - R_{11}D_{21}$ | 0 | 0 | 0 | 0 | 1 | $-1$ |

This matrix encoding of the constraints suggests a similar encoding for the parametrization.

NOTATION 3.7. Let $\mathcal{A}^{(n_1,\ldots,n_m)}$ be a matrix that has one column for each ordered pair $(u, v) \in E \times E$ with $u \sim v$ and one row for each of the parameters identified in Remark 3.3; the column corresponding to the pair $(u, v)$ has $+1$ in the coordinate corresponding to the parameter $a_u$, $-1$ in the coordinate corresponding to the parameter $a_v$ and $+1$ in the coordinate corresponding to the parameter $W(u, v)$. Such a column records the parametrization of $P(u, v)$ in (3.2).

REMARK 3.8. The matrix $\mathcal{A}^{(n_1,\ldots,n_m)}$ has format

$$\left( \prod_{i=1}^{m}(n_i + 1) + \sum_{i=1}^{m} n_i \right) \times \left( 2 \sum_{i=1}^{m} n_i \prod_{j \neq i}(n_j + 1) \right).$$

EXAMPLE 3.9 ($2 \times 1$ grid parametrization matrix). The $9 \times 14$-matrix $\mathcal{A}^{(2,1)}$ has rank 8 and equals

|        | $R_{00}$ | $R_{01}$ | $R_{10}$ | $R_{11}$ | $L_{10}$ | $L_{11}$ | $L_{20}$ | $L_{21}$ | $U_{00}$ | $U_{10}$ | $U_{20}$ | $D_{01}$ | $D_{11}$ | $D_{21}$ |
|--------|------|------|------|------|------|------|------|------|------|------|------|------|------|------|
| $a_{00}$ | 1  | 0  | 0  | 0  | −1 | 0  | 0  | 0  | 1  | 0  | 0  | −1 | 0  | 0  |
| $a_{01}$ | 0  | 1  | 0  | 0  | 0  | −1 | 0  | 0  | −1 | 0  | 0  | 1  | 0  | 0  |
| $a_{10}$ | −1 | 0  | 1  | 0  | 1  | 0  | −1 | 0  | 0  | 1  | 0  | 0  | −1 | 0  |
| $a_{11}$ | 0  | −1 | 0  | 1  | 0  | 1  | 0  | −1 | 0  | −1 | 0  | 0  | 1  | 0  |
| $a_{20}$ | 0  | 0  | −1 | 0  | 0  | 0  | 1  | 0  | 0  | 0  | 1  | 0  | 0  | −1 |
| $a_{21}$ | 0  | 0  | 0  | −1 | 0  | 0  | 0  | 1  | 0  | 0  | −1 | 0  | 0  | 1  |
| $h_1$    | 1  | 1  | 0  | 0  | 1  | 1  | 0  | 0  | 0  | 0  | 0  | 0  | 0  | 0  |
| $h_2$    | 0  | 0  | 1  | 1  | 0  | 0  | 1  | 1  | 0  | 0  | 0  | 0  | 0  | 0  |
| $v_1$    | 0  | 0  | 0  | 0  | 0  | 0  | 0  | 0  | 1  | 1  | 1  | 1  | 1  | 1  |

The following result is essentially a restatement of Theorem 3.1. It will enable us to compute the maximal number of linearly independent rows of $\mathcal{S}^{(n_1,\ldots,n_m)}$ from the rank of $\mathcal{A}^{(n_1,\ldots,n_m)}$.

COROLLARY 3.10. *The two vector spaces spanned by the rows of the matrix $\mathcal{S}^{(n_1,\ldots,n_m)}$ and the rows of the matrix $\mathcal{A}^{(n_1,\ldots,n_m)}$ are the orthogonal complements of each other.*

PROOF. First note that each row vector of $\mathcal{A}$ is orthogonal to the row space of $\mathcal{S}$. This follows directly from the parametrization given in Theorem 3.1, so it suffices to prove that any vector $w \in \ker(\mathcal{S})$ is in the rowspan of $\mathcal{A}$. Let $k \times l$ denote the size of the matrix $\mathcal{A}$. Take $w \in \ker(\mathcal{S})$ and consider $g = (e^{w_i})_{i=1,\ldots,l}$, the componentwise exponentiation of $w$. We will construct a transition matrix $P$ from $w$. The transition matrix $P$ has the nonzero entries $P(u,v)$ $u \sim v$ equal to the corresponding entries of $g$. Since $w \in \ker(\mathcal{S})$, the matrix $P$ can be decomposed into a sum of transition matrices $P_1,\ldots,P_m$ (one in each coordinate direction) which commute pairwise. So by Theorem 3.1, we have $g = (u^{a_1},\ldots,u^{a_l})$ where $u$ is the vector of vertex and edge weight parameters, $a_1,\ldots,a_l$, are the columns of $\mathcal{A}$, and we take componentwise exponentiation. We conclude that $w = \log(g) = (\log(t_1),\ldots,\log(t_k))\mathcal{A}$ is a linear combination of the rows of $\mathcal{A}$. □

LEMMA 3.11. *The rank of the matrix $\mathcal{A}^{(n_1,\ldots,n_m)}$ is one less than the number of rows, that is,*

$$(3.5) \qquad \prod_{i=1}^{m}(n_i+1) + \sum_{i=1}^{m} n_i - 1.$$

Proof. By [12], Remark 0.4.6(d), it suffices to show that the Gram matrix $\mathcal{A}^{(n_1,\ldots,n_m)}\mathcal{A}^{(n_1,\ldots,n_m)\top}$ has rank equal to the number in (3.5). Suppose that $x$ is a row of $\mathcal{A}^{(n_1,\ldots,n_m)}$ corresponding to the parameter $a_u$ for a vertex $u$ of the grid $E$. The inner product $x \cdot x$ is twice the vertex degree, that is, the number of $v \in E$ with $u \sim v$—because $a_u$ appears in the two transition probabilities $P(u,v)$ and $P(v,u)$ for each such $v$. If $x'$ and $x''$ are two such rows, then the inner product $x' \cdot x''$ is 0 if $u \not\sim v$ and it is $-2$ if $u \sim v$.

Now suppose that $y$ is a row of $\mathcal{A}^{(n_1,\ldots,n_m)}$ corresponding to one of the edge parameters $W(he_k,(h+1)e_k)$. The inner product $y \cdot y$ is $2\prod_{\ell \neq k}(n_\ell + 1)$ because the edge parameter appears in the transition probabilities $P(he_k + r,(h+1)e_k + r)$ and $P((h+1)e_k + r, he_k + r)$ for $\prod_{\ell \neq k}(n_\ell + 1)$ choices of $r \in \sum_{\ell \neq k}\mathbb{Z}e_\ell$. If $y'$ and $y''$ are two such rows, then $y' \cdot y'' = 0$, since each transition probability only involves a single edge parameter. For the same reason, $y \cdot z = 0$ if $z$ is a row of $\mathcal{A}^{(n_1,\ldots,n_m)}$ corresponding to an edge parameter $W(je_\ell,(j+1)e_\ell)$ with $\ell \neq k$.

Lastly, if $x$ and $y$ are as above, then $x \cdot y$ is clearly 0 unless the vertex $u$ associated with $x$ is in the set $\{he_k + r : r \in \sum_{\ell \neq k}\mathbb{Z}e_\ell\} \cup \{(h+1)e_k + r : r \in \sum_{\ell \neq k}\mathbb{Z}e_\ell\}$ of starting and ending points of transitions that involve the edge parameter associated with $y$. If the vertex $u$ is in this set, then $x \cdot y = (+1 \times +1) + (+1 \times -1) = 0$ also, with the $+1$ coming from a transition probability of the form $P(u,v)$, and the $-1$ coming from a transition probability of the form $P(v,u)$.

Let $\mathcal{X}$ and $\mathcal{Y}_1,\ldots,\mathcal{Y}_m$ be the submatrices of $\mathcal{A}^{(n_1,\ldots,n_m)}$ that correspond to the rows for the vertices of the grid graph $E$, and the parallel edges in direction $e_1,\ldots,e_m$, respectively. We have shown that $\mathcal{Y}_i\mathcal{Y}_i^\top = 2\prod_{j \neq i}(n_j + 1)I$ where $I$ is the identity matrix of appropriate size, and $\mathcal{X}\mathcal{X}^\top = 2(\Delta - J)$ where $\Delta$ is the diagonal matrix listing the vertex degrees and $J$ is the adjacency matrix of the graph. Thus, $\mathcal{X}\mathcal{X}^\top$ is a multiple of the *Laplacian* matrix of the graph, and because the graph is connected, its rank is one less than the number of vertices [9], Lemma 13.1.1. The matrices $\mathcal{X}\mathcal{Y}_i^\top$ and $\mathcal{Y}_i\mathcal{Y}_j^\top$ are all 0 when $i \neq j$. We conclude that $\mathcal{A}^{(n_1,\ldots,n_m)}\mathcal{A}^{(n_1,\ldots,n_m)\top}$ is a block diagonal matrix, with a block of rank $\prod_{i=1}^{m}(n_i + 1) - 1$, one less than its size, and $m$ blocks of full rank $n_1,\ldots,n_m$, respectively. Therefore, the rank of $\mathcal{A}^{(n_1,\ldots,n_m)}\mathcal{A}^{(n_1,\ldots,n_m)\top}$ is precisely the quantity in (3.5).   □

Corollary 3.10 and Lemma 3.11 imply the following result.

Proposition 3.12. *The rank of the constraint matrix* $\mathcal{S}^{(n_1,\ldots,n_m)}$ *equals*

$$(3.6) \qquad 2\sum_{i=1}^{m}\left(n_i\prod_{j\neq i}(n_j+1)\right) - \prod_{i=1}^{m}(n_i+1) - \sum_{i=1}^{m}n_i + 1.$$

REMARK 3.13. The number (3.6) has the following geometric interpretation. We regard $E$ as the edge graph of a cubical cell complex. The $\ell$-cells of that complex have the form

$$\{u + \lambda_1 e_{k_1} + \lambda_2 e_{k_2} + \cdots + \lambda_\ell e_{k_\ell} : 0 \leq \lambda_1, \lambda_2, \ldots, \lambda_\ell \leq 1\}$$

for $u \in E$ and basis vectors $e_{k_1}, \ldots, e_{k_\ell}$ such that $u + \sum_{j=1}^{\ell} e_{k_j} \in E$. It can be shown that

$$(3.6) = \sum_{\ell=2}^{m} (-1)^\ell (\ell+1) \# \ell\text{-cells},$$

but we omit the proof. The reasoning behind this alternating formula, which is reminiscent of an Euler characteristic, can be used to select a row basis of $\mathcal{S}^{(n_1, \ldots, n_m)}$.

**4. The toric ideal.** In this section we examine the matrix $\mathcal{A} = \mathcal{A}^{(n_1, \ldots, n_m)}$ from the perspective of combinatorial commutative algebra [6, 20]. Let $\mathbb{R}[P]$ be the polynomial ring over the real numbers $\mathbb{R}$ generated by the unknowns $P(u, v)$ with $u, v \in E$ and $u \sim v$. We write $I_{\mathcal{A}}$ for the *toric ideal* associated with the matrix $\mathcal{A}$. Thus $I_{\mathcal{A}}$ is generated by all binomials $P^{Z_+} - P^{Z_-}$ where we take componentwise exponentiation, and $Z = Z_+ - Z_-$ runs over all integer vectors in the kernel of $\mathcal{A}$. Among these binomials are the quadratic binomials in (3.1) whose corresponding vectors $Z$ are the rows of the matrix $\mathcal{S} = \mathcal{S}^{(n_1, \ldots, n_m)}$. These quadratic binomials do not suffice to generate the toric ideal $I_{\mathcal{A}}$, and one of our objectives is to identify a generating set of binomials. In algebraic statistics [5], Section 1.3, one considers several generating sets of $I_{\mathcal{A}}$. Minimal generating sets are known as *Markov bases*. However, it is often more natural to study the larger *Graver basis*, which also contains all circuits of the integer kernel of $\mathcal{A}$, and all reduced Gröbner bases of $I_{\mathcal{A}}$. Recall that the *circuits* of a sublattice of $\mathbb{Z}^N$ are the primitive vectors of inclusion-minimal support. Here "primitive" means that the coordinates have greatest common denominator equal to one, so there are finitely many circuits, and they are unique up to sign. For an elementary introduction to toric ideals we refer to [20], and for their interpretation in terms of Markov chains, see [3, 5].

EXAMPLE 4.1 ($2 \times 1$ grid). If $m = 2$, $n_1 = 2$ and $n_2 = 1$, then $\mathbb{R}[P]$ is a polynomial ring in 14 unknowns and $I_{\mathcal{A}}$ is the ideal of $2 \times 2$-minors of the two matrices in (2.2). Thus the Markov basis of $\mathcal{A}$ has 12 elements. The Graver basis of $\mathcal{A}$ has 29 elements, listed in Example 4.12 below.

The space of commuting birth-and-death processes with positive probabilities is the set of positive real points on a toric variety. The ideal representing that toric variety is the toric ideal $I_{\mathcal{A}}$. The toric variety $V(I_{\mathcal{A}})$ is $D$-dimensional in $\mathbb{R}^N$ where $N = 2\sum_{i=1}^{m} n_i \prod_{j \neq i}(n_j + 1)$ and $D = \prod_{i=1}^{m}(n_i + 1) + \sum_{i=1}^{m} n_i - 1$. These numbers were derived in Section 3. Moreover, a parametrization of the toric variety of commuting birth-and-death processes was given in (3.2).

We start our study in this section by reviewing two definitions regarding matrices in general.

DEFINITION 4.2. A *Cayley matrix* is a special block matrix of the form

$$\begin{bmatrix} \boxed{*} & \boxed{*} & \boxed{*} & \boxed{*} \\ 1 \cdots 1 & 0 & 0 & 0 \\ 0 & 1 \cdots 1 & 0 & 0 \\ 0 & 0 & \ddots & 0 \\ 0 & 0 & 0 & 1 \cdots 1 \end{bmatrix}.$$

By changing the order of the columns of our matrix $\mathcal{A}$, we can see that $\mathcal{A}$ is a Cayley matrix where the blocks correspond to the parallel edges in the lattice. Toric varieties associated with Cayley matrices have a nice, special structure, and they appear in applied contexts such as chemical reaction networks [2], Theorem 9. Here is an even more important special structure.

DEFINITION 4.3. An integer matrix of rank $d$ is *unimodular* if all its nonzero $d \times d$ minors have the same absolute value. There are many equivalent characterizations. For instance, a matrix is unimodular if and only if all initial monomial ideals of its toric ideal are squarefree [20], Section 10.

Investigating whether a matrix is unimodular is interesting from the perspective of algebraic statistics and computational algebra. One reason is the following result which is proven in [20].

PROPOSITION 4.4. *Let $A$ be any unimodular matrix. Then every reduced Gröbner bases of the toric ideal $I_A$ consists of differences of squarefree monomials. Moreover, the following three sets coincide: the union of all reduced Gröbner bases, the set of circuits and the Graver basis of $A$.*

It is thus natural to ask whether our matrix is unimodular. We address this question as follows:

THEOREM 4.5. *The Cayley matrix $\mathcal{A} = \mathcal{A}^{(n_1,\ldots,n_m)}$ is unimodular if and only if the dimension of the grid equals $m = 2$ and the format of the grid is either $2 \times 2$ or $n \times 1$ for some $n \geq 1$.*

Our proof of this theorem employs tools from matroid theory. Recall that a *matroid* $M$ is a pair $(\mathcal{E}, \mathcal{I})$ consisting of a finite ground set $\mathcal{E}$ and a collection $\mathcal{I}$ of subsets of $\mathcal{E}$ that satisfy:

(i) $\varnothing \in \mathcal{I}$.

(ii) If $I_1 \in \mathcal{I}$ and $I_2 \subseteq I_1$, then $I_2 \in \mathcal{I}$.

(iii) If $I_1, I_2 \in \mathcal{I}$ and $\#I_1 < \#I_2$; then there is an element $e \in I_2 \setminus I_1$ such that $I_1 \cup \{e\} \in \mathcal{I}$.

The members of $\mathcal{I}$ are the *independent sets* of $M$, and subsets of $\mathcal{E}$ that are not in $\mathcal{I}$ are *dependent*. A maximal independent set in $M$ is a *basis* of $M$, and a minimal dependent set a *circuit* of $M$. A particular class of matroids arises from matrices as follows. The ground set $\mathcal{E}$ consists of the columns of a matrix $A$, and $\mathcal{I}$ is the set of linearly independent subsets of column vectors of $A$.

The following three results from matroid theory will be useful for the proof of Theorem 4.5. These and further properties of matroids can be found in [18].

LEMMA 4.6. *Let $\mathcal{E}$ be a set and $\mathcal{C}$ a collection of subsets of $\mathcal{E}$. Then $\mathcal{C}$ is the collection of circuits of a matroid on $\mathcal{E}$ if and only if $\mathcal{C}$ satisfies the following three conditions:*

(i) $\varnothing \notin \mathcal{C}$.

(ii) *If $C_1, C_2 \in \mathcal{C}$ and $C_1 \subseteq C_2$, then $C_1 = C_2$.*

(iii) *If $C_1, C_2 \in \mathcal{C}$ distinct and $e \in C_1 \cap C_2$, then $\exists C_3 \in \mathcal{C}$ such that $C_3 \subseteq (C_1 \cup C_2) \setminus \{e\}$.*

LEMMA 4.7. *The collection $\mathcal{B}$ of bases of a matroid $M$ satisfies the following two conditions:*

(i) $\mathcal{B}$ *is nonempty.*

(ii) *If $B_1, B_2 \in \mathcal{B}$ and $e \in B_1 \setminus B_2$, then $\exists f \in B_2 \setminus B_1$ such that $(B_1 \setminus \{e\}) \cup \{f\} \in \mathcal{B}$.*

LEMMA 4.8. *Let $M$ be a matroid over a set $\mathcal{E}$ and $B$ a basis of $M$. If $e \in \mathcal{E} \setminus B$, then $B \cup \{e\}$ contains a unique circuit $C(e, B)$. Moreover, $e \in C(e, B)$.*

Our goal is to describe the circuits of the matroid associated with the matrix $\mathcal{A} = \mathcal{A}^{(n_1, \ldots, n_m)}$. The ground set $\mathcal{E}$ is the set of column vectors of $\mathcal{A}$. We identify $\mathcal{E}$ with the set of directed edges in the grid $E$. Two directed edges in $\mathcal{E}$ are called *parallel* if they point in the same direction, and they have the same edge parameter $W(u, v)$; that is, directed edges $(u', v')$ and $(u'', v'')$ are parallel if $v' - u' \in \{\pm e_k\}$ for some $k$, and $(u'', v'') = (u' + w, v' + w)$ for some $w \in \sum_{\ell \neq k} \mathbb{Z} e_\ell$. Recall that a *cycle* in a graph that is not necessarily simple is a subgraph in which every vertex has degree 2.

DEFINITION 4.9. Let $\mathcal{C} = \mathcal{C}^{(n_1, \ldots, n_m)}$ be the set of subsets $C \subset \mathcal{E}$ with the following properties:

(i) $C$ is a disjoint union of pairs of distinct parallel directed edges.

(ii) The set of undirected edges corresponding to $C$, where two directed edges $(u, v)$ and $(v, u)$ are replaced by two undirected edges, is a union of edge disjoint cycles.

(iii) If $B \subset C$, then at least one of (i) or (ii) does not hold for $B$.

The set $\mathcal{C}^{(2,1)}$ of such walks on the $2 \times 1$-grid has 29 elements, listed in Example 4.12 below.

PROPOSITION 4.10. *The set $\mathcal{C}$ is the set of circuits of a matroid on $\mathcal{E}$.*

PROOF. The collection $\mathcal{C}$ satisfies the three conditions in Lemma 4.6: Condition (i) is trivial, (ii) comes from the minimality requirement (iii) in the definition of $\mathcal{C}$, and (iii) holds because by eliminating the shared directed edges in two overlapping elements, $C_1, C_2 \in \mathcal{C}$, we get a collection of directed edges, which satisfies (i) and (ii) of Definition 4.9 and hence contains an element of $\mathcal{C}$. Thus $\mathcal{C}$ is the set of circuits of a matroid on $\mathcal{E}$.  □

In an early stage of our project, we believed that Proposition 4.10 characterizes the matroid of $\mathcal{A}$. Later we found out that this is true over a field of characteristic two but generally false over the real field $\mathbb{R}$ we are interested in here. However, when $\mathcal{A}$ is unimodular the matroid is independent of the characteristic of the ground field. Here is one instance where this happens.

LEMMA 4.11. *For a grid of size $n \times 1$, the set $\mathcal{C}$ equals the set of circuits of the matroid of $\mathcal{A}$.*

PROOF. For $n \times 1$ grids, the condition (iii) in Definition 4.9 implies that the set of undirected edges corresponding to a circuit $C \in \mathcal{C}$ is a union of, at most, two cycles. For each $C \in \mathcal{C}$ we construct a vector $V_C$ in $\{0, +1, -1\}^{\mathcal{E}}$ whose support equals $C$ and that lies in the kernel of $\mathcal{A}$. The vector $V_C$ is constructed as follows: we choose a cycle in $C$ and walk along this cycle clockwise. We assign $+1$ to forward-pointing edges and $-1$ to backward-pointing edges. If $C$ consists of two cycles, we walk along the second cycle counterclockwise and assign $+1$ to forward-pointing edges and $-1$ to backward-pointing edges. This ensures that $V_C$ has zero inner product with each row of $\mathcal{A}$ that is indexed by a vertex parameter $a_v$. Parallel pairs of edges receive opposite signs because they point in the same direction. This ensures that $V_C$ has zero inner product with each row of $\mathcal{A}$ that is indexed by a parameter $W(u, v)$. Hence $V_C$ lies in the kernel of $\mathcal{A}$.

The construction of $V_C$ reveals that the columns of $\mathcal{A}$ corresponding to any proper subset of $C$ are linearly independent. This implies the inclusion

$$\mathcal{C} \subseteq \{\text{minimal linearly dependent subsets of column vectors of } \mathcal{A}\}.$$

But the reverse inclusion $\supseteq$ holds as well. Consider any $V \in \ker(\mathcal{A}) \setminus \{0\}$. We interpret $V$ as a multiset of signed directed edges in $\mathcal{E}$ where each entry denotes the number of the corresponding directed edges. These edges come in parallel pairs of opposite sign because $V$ has zero inner product with the rows of $\mathcal{A}$ that are indexed by parameters $W(u, v)$. The undirected graph underlying $V$ is a union of cycles because $V$ has zero inner product with the other rows. Hence the support of $V$ must contain some circuit $C \in \mathcal{C}$. This completes the proof of Lemma 4.11.  □

PROOF OF THEOREM 4.5.  To prove the "only if" direction, it suffices to give two examples of circuits which are not squarefree, one for the $3 \times 2$ grid and one for the $1 \times 1 \times 1$ grid. These examples are also circuits for larger two-dimensional grids, and for grids of dimension $m \geq 3$, respectively. Then, by Proposition 4.4, we conclude that the corresponding matrix $\mathcal{A}$ is not unimodular.

(i) $3 \times 2$ grid: the following nonsquarefree binomial is a circuit:

$$L_{10} L_{30} U_{01}^2 U_{20} D_{02} D_{31} - L_{12} L_{32} U_{10} U_{11} U_{31} D_{01} D_{22}.$$

(ii) $1 \times 1 \times 1$ grid: the nonsquarefree binomial

$$R_{000}^2 D_{010} F_{100} B_{001} - R_{001} R_{011} D_{110} F_{010} B_{111}$$

is a circuit where $B$ and $F$ denote the additional backward and forward moves.

We now prove the "if" direction for the $n \times 1$ grid. Let $k$ denote the number of rows of the matrix $\mathcal{A} = \mathcal{A}^{(n,1)}$. Then $\text{rk}(\mathcal{A}) = k - 1$ by Lemma 3.11. We must show that all nonzero $(k-1) \times (k-1)$ minors of $\mathcal{A}$ have the same absolute value. The row vectors of $\mathcal{A}$ corresponding to the vertices of the grid are not independent; they sum to 0. So in order to get nonzero minors, we delete the row corresponding to the vertex $a_{0,\dots,0}$. The resulting matrix is denoted $\mathcal{A}'$.

In Lemma 4.11 we characterized the circuits of the matroid on the column vectors of $\mathcal{A}'$ for $n \times 1$ grids. Consider two distinct bases $B$ and $B'$ of this matroid. To finish the proof, we need to show that the determinants of the submatrices of $\mathcal{A}'$ corresponding to these bases have the same absolute value.

Let $e \in B \setminus B'$. So by Lemma 4.7, $\exists f \in B' \setminus B$ such that $B_1 = (B \setminus \{e\}) \cup \{f\}$ is a basis. This construction can be carried on until $B_k = B'$. So without loss of generality, $B_1 = B'$. By Lemma 4.8, $B \cup \{f\}$ contains a unique circuit $C$ and $f \in C$. Because $B_1$ is also a basis, we get by the same argument that $e \in C$. The corresponding vector $V_C$ in kernel$(\mathcal{A})$ has coordinates $\pm 1$ in positions $e$ and $f$. Hence we can replace the column vector $e$ by the column vector $f$ without changing the absolute value of the determinant. In symbols, $|\det(B)| = |\det(B')|$.

To finish the proof, we need to show that the $13 \times 24$ matrix $\mathcal{A}$ for the $2 \times 2$ grid is unimodular. This can be verified with the help of a computer algebra system such as `CoCoA` or `4ti2`. Namely, we compute the Graver bais of $\mathcal{A}$, and we observe that all binomials in that Graver basis are squarefree. Since the Graver basis contains every reduced Gröbner basis, this implies that every initial monomial of $I_{\mathcal{A}}$ is squarefree. Hence, by [20], Remark 8.10, the matrix $\mathcal{A}^{(2,2)}$ is unimodular. $\square$

We conducted a further computational study of the toric ideal $I_{\mathcal{A}}$ for the two minimal grids whose matrices $\mathcal{A}$ are not unimodular, namely the $3 \times 2$ grid and the $1 \times 1 \times 1$ grid. In what follows, we present two tables that show the total number of binomials of each degree for the Graver basis, the set of circuits and the minimal Markov basis. In parenthesis are the numbers of squarefree binomials in these bases. The two tables were computed using the software `4ti2`.

$3 \times 2$ grid:

| | Degree | | | | | | |
|---|---|---|---|---|---|---|---|
| | 2 | 3 | 4 | 5 | 6 | 7 | 8 |
| Graver basis | 45 (45) | 128 (128) | 464 (464) | 1600 (1600) | 3904 (3904) | 4928 (4032) | 1088 (192) |
| Circuits | 45 (45) | 128 (128) | 464 (464) | 1600 (1600) | 3904 (3904) | 4928 (4032) | 896 (0) |
| Markov basis | 36 (36) | 0 (0) | 4 (4) | 4 (4) | 4 (4) | 0 (0) | 0 (0) |

$1 \times 1 \times 1$ grid:

| | Degree | | | | |
|---|---|---|---|---|---|
| | 2 | 3 | 4 | 5 | 6 |
| Graver basis | 42 (42) | 224 (224) | 1032 (1032) | 1728 (1152) | 672 (96) |
| Circuits | 42 (42) | 224 (224) | 1032 (1032) | 1728 (1152) | 576 (0) |
| Markov basis | 33 (33) | 8 (8) | 12 (12) | 0 (0) | 0 (0) |

Note that the Markov basis in both examples consists of squarefree binomials only and that the number of binomials with squares of degree 7 and 8 for the $2 \times 3$ grid, respectively of degree 5 and 6 for the $1 \times 1 \times 1$ grid is the same. We wonder if this also holds for higher dimensions.

We end this section by illustrating the conclusion of Theorem 4.5 for our running example.

EXAMPLE 4.12 ($2 \times 1$ grid). The set $\mathcal{C}^{(2,1)}$ of circuits of the $8 \times 14$-matrix $\mathcal{A}^{(2,1)}$, in Example 3.9, consists of 29 elements. Each circuit $C$ corresponds to a unique (up to sign) vector $V_C$ in $\{0, +1, -1\}^{14} \cap \text{kernel}(\mathcal{A}^{(2,1)})$, and hence to a difference of squarefree monomials in

$$\mathbb{R}[P] = \mathbb{R}[R_{00}, R_{01}, R_{10}, R_{11}, L_{10}, L_{11}, L_{20}, L_{21}, U_{00}, U_{10}, U_{20}, U_{01}, D_{11}, D_{21}].$$

We now list the 29 binomials in the Graver basis $\mathcal{C}^{(2,1)}$ of $I_{\mathcal{A}^{(2,1)}}$. First, there are the eight quadrics which generate $I^{(2,1)}$. These correspond to the eight circuits of type $C_1$. Next, there are five circuits of type $C_2$. These are displayed below as a quadratic binomial and as a picture:

(1) $R_{00}L_{10} - R_{01}L_{11}$  ,

(2) $R_{10}L_{20} - R_{11}L_{21}$  ,

(3) $D_{11}U_{10} - D_{01}U_{00}$  ,

(4) $D_{21}U_{20} - D_{01}U_{00}$  ,

(5) $D_{21}U_{20} - D_{11}U_{10}$  .

The remaining 16 circuits are given below with the corresponding binomials in the Graver basis.

(1) $R_{00}R_{10}U_{20} - R_{01}R_{11}U_{00}$  ,

(2) $L_{10}L_{20}U_{00} - L_{11}L_{21}U_{20}$  ,

(3) $R_{00}L_{21}U_{20} - R_{01}L_{20}U_{00}$ 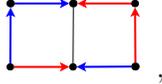 ,

(4) $R_{10}L_{11}U_{20} - R_{11}L_{10}U_{00}$ 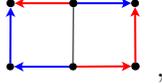 ,

(5) $R_{10}D_{11}U_{20} - R_{11}D_{01}U_{00}$ 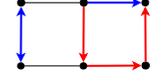 ,

(6) $L_{21}D_{11}U_{20} - L_{20}D_{01}U_{00}$ 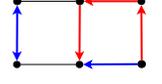 ,

(7) $R_{00}D_{21}U_{20} - R_{01}D_{11}U_{00}$ 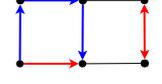 ,

(8) $L_{11}D_{21}U_{20} - L_{10}D_{11}U_{00}$ 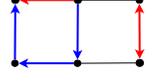 ,

(9) $R_{00}R_{10}D_{01} - R_{01}R_{11}D_{21}$ 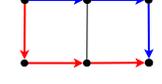 ,

(10) $L_{10}L_{20}D_{21} - L_{11}L_{21}D_{01}$ 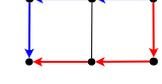 ,

(11) $R_{00}L_{21}D_{01} - R_{01}L_{20}D_{21}$ 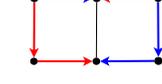 ,

(12) $R_{10}L_{11}D_{01} - R_{11}L_{10}D_{21}$ 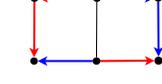 ,

(13) $R_{11}D_{21}U_{10} - R_{10}D_{01}U_{00}$ 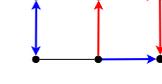 ,

(14) $L_{20}D_{21}U_{10} - L_{21}D_{01}U_{00}$ 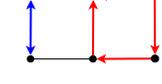 ,

(15) $R_{01}D_{21}U_{20} - R_{00}D_{01}U_{10}$ 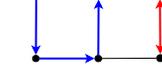 ,

$(16)\ L_{10}D_{21}U_{20} - L_{11}D_{01}U_{10}$  .

Similar lists of circuits can be computed, at least in principle, for larger and higher-dimensional grids. However, their combinatorial structure will be quite complicated. To appreciate this, we invite the reader to draw the picture for the circuit of the $4 \times 4$-grid displayed in (5.6) below.

**5. Primary decomposition and boundary components.** The set of commuting birth-and-death processes is the set of nonnegative real solutions of the quadratic equations (3.1). In the previous two sections we have focused on the set of strictly positive solutions, and we have argued that its closure is a toric variety $V(I_{\mathcal{A}})$ with interesting combinatorial properties. In this section we apply binomial primary decomposition [3, 6] to classify those nonnegative solutions of (3.1) that are not in the closure of the positive solutions.

NOTATION 5.1. Let $I^{(n_1,\ldots,n_m)}$ denote the ideal in the polynomial ring $\mathbb{K}[P]$ that is generated by the quadrics in (3.1). Here $\mathbb{K}$ is allowed to be an arbitrary field. The ideal $I^{(n_1,\ldots,n_m)}$ is a binomial ideal in the sense of [6]. It represents all commuting birth-and-death processes.

To keep things simple, we first concentrate on the two-dimensional case. We use the notation of Section 2 for an $m \times n$ grid. Then $I^{(m,n)}$ is a binomial ideal in the polynomial ring $\mathbb{K}[R, L, D, U]$ over an arbitrary field $\mathbb{K}$. The set of indeterminates equals

$$\{R_{ij} : 0 \leq i < m, 0 \leq j \leq n\} \cup \{L_{ij} : 0 < i \leq m, 0 \leq j \leq n\}$$
$$\cup \{D_{ij} : 0 \leq i \leq m, 0 < j \leq n\} \cup \{U_{ij} : 0 \leq i \leq m, 0 \leq j < n\}.$$

The ideal $I^{(m,n)}$ is minimally generated by the following set of $4mn$ quadratic binomials

$$\{U_{i,j}R_{i,j+1} - R_{i,j}U_{i+1,j}, R_{i,j+1}D_{i+1,j+1} - D_{i,j+1}R_{i,j},$$
$$D_{i+1,j+1}L_{i+1,j} - L_{i+1,j+1}D_{i,j+1},$$
$$L_{i+1,j}U_{i,j} - U_{i+1,j}L_{i+1,j+1} : 0 \leq i < m, 0 \leq j < n\}.$$

As a warm-up, let us discuss the simplest case, namely the four binomials on a single square.

EXAMPLE 5.2 $[(m = n = 1)]$.   The variety $V(I^{(1,1)}) \subset \mathbb{K}^8$ consists of all pairs of matrices

$$(5.1) \quad \begin{pmatrix} 0 & 0 & R_{00} & 0 \\ 0 & 0 & 0 & R_{01} \\ L_{10} & 0 & 0 & 0 \\ 0 & L_{11} & 0 & 0 \end{pmatrix} \quad \text{and} \quad \begin{pmatrix} 0 & U_{00} & 0 & 0 \\ D_{01} & 0 & 0 & 0 \\ 0 & 0 & 0 & U_{10} \\ 0 & 0 & D_{11} & 0 \end{pmatrix}$$

that commute. The possibilities for this to happen are revealed by the primary decomposition,

$$\begin{aligned} I^{(1,1)} = \langle &U_{00}R_{01} - R_{00}U_{10}, R_{01}D_{11} \\ & - D_{01}R_{00}, D_{11}L_{10} - L_{11}D_{01}, L_{10}U_{00} - U_{10}L_{11} \rangle \\ = &I_{\mathcal{A}} \cap \langle U_{00}, U_{10}, D_{01}, D_{11} \rangle \cap \langle R_{00}, R_{01}, L_{10}, L_{11} \rangle. \end{aligned}$$

Here $I_{\mathcal{A}}$ is the toric ideal of the matrix $\mathcal{A} = \mathcal{A}^{(1,1)}$. It is generated by the six $2 \times 2$-minors of the matrix,

$$(5.2) \quad \begin{pmatrix} R_{00} & U_{00} & L_{11} & D_{11} \\ R_{01} & U_{10} & L_{10} & D_{01} \end{pmatrix}.$$

We see that $V(I^{(1,1)})$ is the union of three irreducible components. The main component $V(I_{\mathcal{A}})$ has dimension 5 while the two boundary components have dimension 4. The two matrices (5.1) commute if and only if either one of them is the zero matrix, or the matrix (5.2) has rank $\leq 1$.

An ideal in a polynomial ring is *radical* if and only if it is a finite intersection of prime ideals. We say that a *pure toric ideal* is a prime ideal which is generated by indeterminates and differences of monomials. For example, $\langle xy - z^2 \rangle$ is pure toric, but $\langle xy + z^2 \rangle$ and $\langle x^2 - z^2 \rangle$ are not.

CONJECTURE 5.3.   *The ideal $I^{(m,n)}$ is radical, all its associated primes are pure toric ideals and the prime decomposition of $I^{(m,n)}$ in $\mathbb{K}[R, L, D, U]$ is independent of the coefficient field $\mathbb{K}$.*

Example 5.2 establishes this conjecture for $m = n = 1$. What follows is devoted to establishing partial results that support Conjecture 5.3 and to explore the combinatorics of the irreducible components of our binomial variety $V(I^{(m,n)})$ of commuting birth-and-death processes.

EXAMPLE 5.4 $[(m \leq 3, n = 1)]$.   The radical ideal $I^{(2,1)}$ is the intersection of 11 pure toric ideals. These 11 primes are easily found in $\mathtt{Singular}$ or $\mathtt{Macaulay2}$. They come in six symmetry classes:

1. The main component has codimension 6 and degree 16. It corresponds to the rank one condition for two matrices as in (5.2), one for each of the two squares in the $2 \times 1$-grid,

$$I_{\mathcal{A}} = \left\langle 2 \times 2\text{-minors of } \begin{pmatrix} R_{00} & U_{00} & L_{11} & U_{11} \\ R_{01} & U_{10} & L_{10} & D_{01} \end{pmatrix} \right.$$

$$\left. \text{and of } \begin{pmatrix} R_{10} & U_{10} & L_{21} & D_{21} \\ R_{11} & U_{20} & L_{20} & D_{11} \end{pmatrix} \right\rangle.$$

2. There are two components of codimension 7 and degree 4, each isomorphic to $V(I^{(1,1)})$ times a four-dimensional coordinate subspace. One of these two primes is the pure toric ideal generated by $\{L_{21}, L_{20}, R_{11}, R_{10}\}$ and the six $2 \times 2$-minors of the matrix (5.2).

3. The monomial component representing horizontal transitions has codimension 6,

$$\langle U_{00}, U_{10}, U_{20}, D_{01}, D_{11}, D_{21} \rangle \qquad \text{} .$$

4. The monomial component representing vertical transitions has codimension 8,

$$\langle R_{00}, R_{01}, R_{10}, R_{11}, L_{10}, L_{11}, L_{20}, L_{21} \rangle \qquad \text{} .$$

5. There is another pair of monomial primes of codimension 8. One of these two ideals is

$$\langle R_{00}, R_{01}, L_{10}, L_{11}, U_{10}, U_{20}, D_{11}, D_{21} \rangle \qquad \text{} .$$

6. The last class consists of four monomial primes of codimension 7, such as

$$\langle R_{10}, L_{21}, U_{00}, U_{10}, D_{01}, D_{11}, D_{21} \rangle \qquad \text{} .$$

A similar computation shows that $I^{(3,1)}$ is also a radical ideal. It is the intersection of 40 pure toric ideals. The list is analogous to that above. Here, the main component $I_{\mathcal{A}}$ is generated by the 18 minors of the three $2 \times 4$-matrices coming from the three squares in the grid.

Current general-purpose implementations of primary decomposition are not able to perform the corresponding computation for the $2 \times 2$ grid. To crack this case, we applied the methods of [3, 6] in an interactive fashion in order to verify Conjecture 5.3 and find the prime decomposition.

EXAMPLE 5.5 [($m = n = 2$)]. The ideal $I^{(2,2)}$ is radical and it is the irredundant intersection of 199 pure toric ideals. The toric ideal $I_{\mathcal{A}}$ representing the main component is minimally generated by 26 binomials. Besides the 24 familiar quadrics which come from the $2 \times 4$-matrices on the four squares of the grid, here we find the following two additional quartic Markov basis elements:

$$
\begin{aligned}
(5.3) \qquad & R_{00}L_{22}U_{20}D_{02} - R_{02}L_{20}U_{00}D_{22} \quad \text{and} \\
& R_{10}L_{12}U_{21}D_{01} - R_{12}L_{10}U_{01}D_{11}.
\end{aligned}
$$

We shall next derive one theoretical result related to Conjecture 5.3. Let $I$ be a *pure binomial ideal* in a polynomial ring $\mathbb{K}[x_1, x_2, \ldots, x_r]$, that is, $I$ is an ideal generated by monomial differences $x_1^{u_1} x_2^{u_2} \cdots x_r^{u_r} - x_1^{v_1} x_2^{v_2} \cdots x_r^{v_r}$. The lattice $\mathcal{L}_I$ associated with $I$ is the sublattice of $\mathbb{Z}^r$ spanned by the vectors $(u_1 - v_1, u_2 - v_2, \ldots, u_r - v_r)$ in any generating set of $I$. For any positive integer $i < r$, the elimination ideal $I \cap K[x_1, \ldots, x_i]$ is also pure binomial, and we have

$$
(5.4) \qquad \mathcal{L}_{I \cap K[x_1, \ldots, x_i]} = \mathcal{L}_I \cap (\mathbb{Z}e_1 + \cdots + \mathbb{Z}e_i).
$$

A result of Gilmer [8] states that the radical $\sqrt{I}$ of a pure binomial ideal $I$ is also a pure binomial ideal, and, using the methods developed in [6], Section 3, it is not hard to see that $\mathcal{L}_I = \mathcal{L}_{\sqrt{I}}$.

We say that a pure binomial ideal $I$ is *unimodular* if the corresponding lattice $\mathcal{L}_I$ is a *unimodular* sublattice of $\mathbb{Z}^r$. By this we mean that the quotient group $\mathbb{Z}^r / \mathcal{L}_I$ is free abelian and the same holds for every elimination ideal of $I$. Equivalently, $I$ is unimodular if the quotient group $\mathbb{Z}\{e_s : s \in \sigma\} / (\mathcal{L}_I \cap \mathbb{Z}\{e_s : s \in \sigma\})$ is free abelian for every subset $\sigma$ of $\{1, 2, \ldots, r\}$.

PROPOSITION 5.6. *If $I$ is a unimodular pure binomial ideal in a polynomial ring $\mathbb{K}[x_1, \ldots, x_r]$, then every associated prime of $I$ is a pure toric ideal.*

PROOF. We abbreviate the polynomial ring by $S = \mathbb{K}[x_1, \ldots, x_r]$. If $I$ is toric then we are done. Otherwise, we follow the approach in the proof of [6], Theorem 6.1, and pick a variable $x_i$ such that both $(I : x_i)$ and $I + \langle x_i \rangle$ strictly contain $I$. Note that the lattice $\mathcal{L}_{I + \langle x_i \rangle}$ equals $\mathcal{L}_I \cap \mathbb{Z}\{e_1, \ldots, e_{i-1}, e_{i+1}, \ldots, e_r\}$ and hence is unimodular. The exact sequence,

$$
(5.5) \qquad 0 \longrightarrow S/(I : x_i) \longrightarrow S/I \longrightarrow S/(I + \langle x_i \rangle) \longrightarrow 0,
$$

shows that each associated prime of $I$ is an associated prime of $(I : x_i)$ or of $I + \langle x_i \rangle$. Both $(I : x_i)$ and $I + \langle x_i \rangle$ are unimodular pure binomial ideals, and they satisfy the desired conclusion by Noetherian induction. Hence so does $I$. $\square$

We have the following result which generalizes the main conclusion of Example 5.4.

COROLLARY 5.7. *Every associated prime of the binomial ideal $I^{(n,1)}$ is a pure toric ideal.*

PROOF. In Theorem 4.5 we have shown that the lattice kernel($\mathcal{A}^{(n,1)}$) associated with our binomial ideal $I^{(n,1)}$ is unimodular. By Proposition 5.6, every associated prime of $I^{(n,1)}$ is pure toric.  □

An overly optimistic conjecture suggested by Examples 5.4 and 5.5 would be that a minimal Markov basis for the main component $I_{\mathcal{A}}$ of $I^{(m,n)}$ consists of the $6mn$ quadrics derived from the matrices (5.2) on all $1 \times 1$-subgrids and the $2(m-1)(n-1)$ quartic binomials (5.3) derived from all $2 \times 2$-subgrids of the $m \times n$-grid. Unfortunately this is far from true. For instance, for $m = n = 3$ the Markov basis of $I_{\mathcal{A}}$ consists of 314 binomials, including the 54 quadrics and 8 quartics, but there are also 16 quintics, 36 sextics and 200 binomials of degree eight, such as

$$(5.6) \qquad U_{00}U_{32}L_{20}L_{30}L_{13}L_{23}R_{21}R_{02} - U_{30}U_{02}L_{11}L_{21}L_{22}L_{32}R_{00}R_{23}.$$

It would be desirable to identify which circuits make up the Markov basis for arbitrary $m$ and $n$.

We next examine the situation for the smallest example of a three-dimensional grid.

EXAMPLE 5.8 ($1 \times 1 \times 1$-grid). The ideal $I^{(1,1,1)}$ is generated by 24 quadratic binomials (four for each face of the 3-cube) in 24 unknowns (two for each edge of the 3-cube). We write

$$\begin{aligned}
I^{(1,1,1)} = \langle\, &L_{100}U_{000} - U_{100}L_{110}, L_{101}U_{001} - U_{101}L_{111}, \\
&L_{110}D_{010} - D_{110}L_{100}, L_{111}D_{011} - D_{111}L_{101}, \\
&L_{100}F_{000} - F_{100}L_{101}, L_{110}F_{010} - F_{110}L_{111}, \\
&L_{101}B_{001} - B_{101}L_{100}, L_{111}B_{011} - B_{111}L_{110}, \\
&R_{000}U_{100} - U_{000}R_{010}, R_{001}U_{101} - U_{001}R_{011}, \\
&R_{010}D_{110} - D_{010}R_{000}, R_{011}D_{111} - D_{011}R_{001}, \\
&R_{000}F_{100} - F_{000}R_{001}, R_{010}F_{110} - F_{010}R_{011}, \\
&R_{001}B_{101} - B_{001}R_{000}, R_{011}B_{111} - B_{011}R_{010}, \\
&B_{101}U_{100} - U_{101}B_{111}, B_{001}U_{000} - U_{001}B_{011}, \\
&B_{011}D_{010} - D_{011}B_{001}, B_{111}D_{110} - D_{111}B_{101},
\end{aligned}$$

$$F_{000}U_{001} - U_{000}F_{010}, F_{100}U_{101} - U_{100}F_{110},$$

$$F_{010}D_{011} - D_{010}F_{000}, F_{110}D_{111} - D_{110}F_{100}\rangle.$$

The main component is a unimodular toric ideal of codimension 14 and degree 300. Its Markov basis consists of 53 binomials (33 quadrics, 8 cubics and 12 quartics), and its Graver basis consists of 3698 binomials (42 quadrics, 224 cubics, 1032 quartics, 1728 quintics and 672 sextics).

We found that the binomial ideal $I^{(1,1,1)}$ is the intersection of 131 prime ideals, so it is radical. Besides the main component, there are 91 monomial primes and 39 other binomial primes:

| # | cd | deg | representative associated prime |
|---|---|---|---|
| 2 | 12 | 1 | $\langle B_{001}, B_{111}, D_{010}, D_{111}, F_{010}, F_{100},$ |
|   |   |   | $L_{100}, L_{111}, R_{001}, R_{010}, U_{001}, U_{100}\rangle$ |
| 24 | 15 | 1 | $\langle B_{011}, B_{101}, B_{111}, D_{011}, D_{110}, F_{000}, F_{100}, F_{110},$ |
|   |   |   | $L_{101}, L_{110}, L_{111}, R_{000}, R_{001}, R_{011}, U_{000}\rangle$ |
| 9 | 16 | 1 | $\langle B_{001}, B_{011}, B_{101}, B_{111}, F_{000}, F_{010}, F_{100}, F_{110},$ |
|   |   |   | $L_{100}, L_{101}, L_{110}, L_{111}, R_{000}, R_{001}, R_{010}, R_{011}\rangle$ |
| 48 | 17 | 1 | $\langle B_{101}, B_{111}, D_{010}, D_{011}, D_{111}, F_{010}, F_{100}, F_{110},$ |
|   |   |   | $L_{100}, L_{101}, L_{110}, L_{111}, R_{000}, R_{001}, R_{010}, U_{000}, U_{001}\rangle$ |
| 8 | 18 | 1 | $\langle B_{011}, B_{101}, B_{111}, D_{010}, D_{011}, D_{110}, F_{010}, F_{100}, F_{110},$ |
|   |   |   | $L_{100}, L_{101}, L_{111}, R_{000}, R_{001}, R_{011}, U_{000}, U_{001}, U_{100}\rangle$ |
| 3 | 14 | 16 | $\langle F_{110}, F_{100}, F_{010}, F_{000}, B_{111}, B_{101}, B_{011}, B_{001}$ |
|   |   |   | and 12 of the binomial generators of $I^{(1,1,1)}\rangle$ |
| 12 | 15 | 4 | $\langle R_{000}, R_{010}, L_{100}, L_{110}, F_{110}, F_{100}, F_{010}, F_{000},$ |
|   |   |   | $B_{111}, B_{101}, B_{011}, B_{001}$ and six binomials$\rangle$ |
| 24 | 16 | 9 | $\langle L_{111}, L_{101}, F_{000}, F_{010}, R_{010}, R_{000}, L_{110}, L_{100},$ |
|   |   |   | $F_{110}, F_{100}, B_{111}, B_{101}$ and six binomials$\rangle$. |

We end by restating the main point of Conjecture 5.3 for grids in the $m$-dimensional lattice.

CONJECTURE 5.9.   *The binomial ideal $I^{(n_1,\dots,n_m)}$ is radical.*

*Note added in proof*: On May 4, 2009, Thomas Kahle announced that he had disproved Conjectures 5.3 and 5.9 using his new software for primary decomposition of binomial ideals. Kahle's result states that the ideal $I^{(2,3)}$ is not radical. More precisely, this ideal is the intersection of 2638 primary ideals of which 11 are not prime. A witness is given by the binomial

$$f = D_{01}R_{03}R_{10}L_{12}U_{21}L_{22}D_{23} - U_{01}R_{03}L_{10}R_{13}D_{21}L_{23}D_{23}.$$

Using `Macaulay 2` or `Singular`, we can check that $f^2 \in I^{(2,3)}$ but $f \notin I^{(2,3)}$.

S. N. Evans
Department of Statistics #3860
University of California at Berkeley
367 Evans Hall
Berkeley, California 94720-3860
USA
E-mail: evans@stat.berkeley.edu
URL: http://www.stat.berkeley.edu/users/evans/

B. Sturmfels
Department of Mathematics #3840
University of California at Berkeley
970 Evans Hall
Berkeley, California 94720-3840
USA
E-mail: bernd@math.berkeley.edu
URL: http://www.math.berkeley.edu/˜bernd/

C. Uhler
Department of Statistics #3860
University of California at Berkeley
345 Evans Hall
Berkeley, California 94720-3860
USA
E-mail: cuhler@stat.berkeley.edu
URL: http://www.carolineuhler.com/